\hsize=120mm
\vsize=185mm
\parindent=8mm
\frenchspacing

\pageno=1 

\def\Bbb #1{{\bf #1}}

\def\remark{{\par\noindent{\bf Remark. }}}
\def\proof{\par\noindent{\bf Proof. }}
\def\sqr#1#2{{\vcenter{\vbox{\hrule height.#2pt
   \hbox{\vrule width.#2pt height#1pt \kern#1pt
     \vrule width.#2pt}
   \hrule height.#2pt}}}}
\def\whis{\mathchoice\sqr34\sqr34\sqr{2.1}3\sqr{1.5}3}
\def\qed{\ifmmode \eqno{\whis}\else
{{\unskip\nobreak\hfil\penalty50\hskip 2em\hbox{}\nobreak\hfil
${\whis}$\parfillskip=0pt\par\medskip}}\fi}
\def\qqed{{\unskip\nobreak\hfil\penalty50\hskip 2em\hbox{}\nobreak\hfil
${\whis\quad\whis}$\parfillskip=0pt\par\medskip}}

\def\ade{1}
\def\add{2}
\def\adc{3}
\def\ce{4}
\def\evf{5}
\def\eve{6}
\def\hen{7}
\def\lea{8}
\def\leb{9}
\def\ven{10}

\centerline{\bf A bound on the exponent of the cohomology of $BC$-bundles}
\centerline{I. J. Leary}
\centerline{Centre de Recerca Matematica,}
\centerline{Institut d'Estudis Catalans,}
\centerline{Apartat 50,}
\centerline{E-08193 Bellaterra.}

\noindent
We give a lower bound for the exponent of certain elements in
the integral cohomology of the total spaces of principal $BC$-bundles for
$C$ a finite cyclic group.  We are mainly interested in the case when 
the total space is $BG$ for some discrete group $G$ having a central subgroup 
isomorphic to $C$.  As applications we give a proof of 
the theorem of A. Adem and H.-W. Henn 
that a $p$-group is elementary abelian if and 
only if its integral cohomology has exponent $p$, and 
we exhibit some infinite groups of finite virtual cohomological 
dimension whose Tate-Farrell cohomology contains torsion of order 
greater than the l.c.m.\ of the orders of their finite subgroups.  
Our examples include a class of groups having similar properties 
discovered by Adem and J. Carlson.  
   As a third application, we 
examine the integral cohomology of a class of $p$-groups expressible 
as central extensions with cyclic kernel and quotient abelian of 
$p$-rank two.  For each such $G$ we determine the minimal $n$ such 
that almost all (i.e.\ all but possibly finitely many) of the 
groups $H\sp i(BG)$ have exponent dividing $p\sp n$.  
The lemma we use to give an upper bound for the exponents of almost
all of the groups $H^i(BG)$ applies to any $p$-group and may be of
independent interest.  
Here, and throughout
the paper, the coefficients for cohomology are to be the integers when 
not otherwise stated, and we write $\Bbb Z_n$ for the integers 
modulo~$n$.  The author gratefully acknowledges that this work was 
funded by the DGICYT.

\proclaim Proposition 1.  Let $C$ be a cyclic group of order $n$, 
and let $E$ be a principal $BC$-bundle over a connected space $X$, 
classified by $\xi\in H^2(X;C)$ of order $m$.  Then for any $i \geq 0$, 
any element of $H^{2i}(E)$ restricting to the fibre as a generator for 
$H^{2i}(BC)$ has order divisible by $mn$.  

\remark Note that we do not claim that such elements always exist, 
nor do we rule out the possibility that they have infinite order.  

\proof In [\ce] Cartan and Eilenberg computed the ring $H^*(BC;R)$ 
for any coefficient ring $R$.  Recall that we have the following 
ring isomorphisms: 
$$H^*(BC)\cong \Bbb Z[z]/(nz),\qquad
H^*(BC;\Bbb Z_{n}) \cong \Bbb Z_{n}[x,y]/(ny,nx,y^2-ex),$$
where $e=0$ if $n$ is odd and $e=n/2$ if $n$ is even, and $y$ has
degree 1 while $x$ and $z$ have degree two.  The natural map from 
integral to mod-$n$ cohomology sends $z$ to $x$, and 
if we let $\beta$ stand for the Bockstein for the coefficient sequence 
$$0\rightarrow \Bbb Z\rightarrow \Bbb Z\rightarrow \Bbb Z_{n}
\rightarrow 0,$$
then it is easy to see that $\beta(y)=z$, and that 
therefore $\beta(yx^i)=z^{i+1}$.  

Now consider the spectral sequence for the given fibration 
with coefficients in $\Bbb Z_{n}$.  By assumption the fundamental group 
of $X$ acts trivially on the cohomology of $BC$, and so 
$$E_2^{i,j}\cong H^i(X;\Bbb Z_{n})\otimes H^j(BC;\Bbb Z_{n}).$$
Now $1\otimes yx^j$ represents a generator for $E_2^{0,2j+1}$ and 
$1\otimes x^j$ represents a generator for $E_2^{0,2j}$.  Comparing 
this spectral sequence with the spectral sequence for the path-loop 
fibration over an Eilenberg-MacLane space $K(C,2)$ it is easy to see 
that $d_2(1\otimes y)= \xi$ and $d_2(1\otimes x)=0$.  (In fact, 
$d_3(1\otimes x) = \xi'\otimes 1$, where $\xi'$ is the image of
$\beta(\xi)$ under the map from $H^3(X)$ to $H^3(X;\Bbb Z_{n})$, and 
$d_4$ may be described using the argument given in [\lea], but 
we do not need this here.)  Now $d_2(1\otimes x^jy)=\xi\otimes x^j$
and $d_2(1\otimes x^j)=0$, from which it follows that $E_3^{0,2j}$ is 
generated by $1\otimes x^j$ and $E_3^{0,2j+1}$ by $m(1\otimes yx^j)$.  
The map from $H^*(E;\Bbb Z_{n})$ to $H^*(BC;\Bbb Z_{n})$ factors through 
$E_\infty^{0,*}$, which is a subgroup of $E_3^{0,*}$, and so 
we see that the image of $H^{2j+1}(E;\Bbb Z_{n})$ in $H^{2j+1}(BC;\Bbb Z_{n})$
must be contained in the subgroup generated by $myx^j$.  

Now recall that the image of the Bockstein $\beta$ defined above is 
exactly the elements of integral cohomology of order dividing $n$.  
Let $f: BC\rightarrow E$ be the inclusion of the fibre of the above 
fibration.  Now let $\chi$ be an element of $H^*(E)$ such that 
$f^*(\chi)=z^{j+1}$ for some $j$.  If $\chi$ has infinite order then 
there is nothing to prove.  Otherwise, the order of $\chi$ must be a 
multiple of $n$ (the order of $z^{j+1}$), say $m'n$, and it remains to 
show that $m$ divides $m'$.  Now $m'\chi$ has order $n$, so there 
exists $\chi'\in H^{2j+1}(E;\Bbb Z_{n})$ such that $\beta(\chi')=m'\chi$.
However, the spectral sequence argument shows that $f^*(\chi')$ is 
in the subgroup of $H^{2j+1}(BC;\Bbb Z_{n})$ generated by $myx^j$ and 
hence $\beta f^*(\chi')$ is in the subgroup of $H^{2j+2}(BC)$ generated 
by $mz^{j+1}$, but $\beta f^*(\chi')= f^*\beta(\chi') 
= f^*(\chi) = m'z^{j+1}$. 
\qed

\proclaim Corollary 1.  Let $C$ be a cyclic subgroup of order~$n$ of a 
group $G$.  If there exists an element of $H^*(BG)$ of order~$n$ whose
image in~$H^*(BC)$ is a generator for $H^{2i}(BC)$ for some~$i$, then 
$C$ is a direct factor of its centraliser in~$G$.  

\proof This is just Proposition~1 applied to the principal $BC$-bundle 
with total space the classifying space of the centraliser of~$C$.  
\qed

\proclaim Corollary~2.  Let $G$ be a discrete group expressible as 
a central extension with kernel~$C$ cyclic of order~$n$.  Let $Q$  be
the quotient~$G/C$, and let the extension class of $G$ in $H^2(BQ;C)$
have order~$m$.  If $G$ has a normal subgroup $N$ of finite index 
whose intersection with~$C$ is trivial (for example, if $G$ is finite 
or residually finite), then for infinitely many~$i$, $H^{2i}(BG)$ 
contains elements of order~$mn$.  

\remark The condition that the extension class of $G$ has order~$m$ 
may be rephrased as follows:  If $D$ is the smallest subgroup of $C$ 
such that $G/D$ is isomorphic to $(C/D)\times Q$, then $D$ has order~$m$.  

\proof Let $G'$ be the quotient $G/N$, and let $C'$ be the image of 
$C$ in $G'$.  Then $C'$ is isomorphic to $C$ and $G'$ is finite.  
By either Evens' argument using the Norm map from $H^*(BC)$ to 
$H^*(BG)$ [\evf,\eve] or Venkov's argument using Chern classes of 
a representation of $G'$ restricting faithfully to $C'$ [\ven],  
we see that for infinitely many $i$ there exists   
$\chi'\in H^{2i}(BG')$ whose image in $H^{2i}(BC')$ is a generator.  
If $\chi$ is the image of $\chi'$ in $H^*(BG)$, then $\chi$ has 
finite order (dividing the order of $G'$) and its image in 
$H^{2i}(BC)$ is a generator.  Hence by Proposition~1, some multiple of 
$\chi$ has order exactly~$mn$.  
\qed

The first example of a group whose Tate-Farrell cohomology contains 
elements of order greater than the l.c.m.\ of the orders of its 
finite subgroups is due to Adem [\add].  The following application of 
Corollary~2 is more closely related to some other examples due to 
Adem and Carlson [\adc].  In particular, Corollary~3 may be compared 
with Theorem~3.1 of [\adc], which gives stronger cohomological information
about a smaller class of groups.  

\def\pra{\par}
\proclaim Corollary 3.  With notation and hypotheses as in Corollary~2, 
assume also that $Q$ has finite cohomological dimension (or equivalently, 
assume that there is a finite-dimensional CW-complex~$BQ$).  Then 
\pra\hangindent=\parindent
a) $G$ has finite virtual cohomological dimension and hence the 
Tate-Farrell cohomology groups $\hat H^i(G)$ are defined,
\pra
b) $C$ consists of all the elements of $G$ of finite order, and 
\pra
c) $\hat H^i(G)$ contains elements of order $mn$ for infinitely many $i$.  

\proof The subgroup $N$ of $G$ has finite index and is isomorphic to 
a subgroup of $Q$, so has cohomological dimension less than or equal 
to that of $Q$.  Hence $G$ has finite vcd.  The group $Q$ is torsion-free,
and so any element of $G-C$ has infinite order because its image in 
$Q$ does.  If $i$ is greater than $\hbox{vcd}G$ then $\hat H^i(G)$ is 
isomorphic to $H^i(BG)$, and so the third claim follows from 
Corollary~2.  
\qed

The following Corollary is due to Adem [\ade] and Henn [\hen].   
\proclaim Corollary 4.  Let $G$ be a finite $p$-group.  Then $G$ is not 
elementary abelian if and only if $H^i(BG)$ contains elements of order~$p^2$ 
for some~$i$ if and only if $H^i(BG)$ contains elements of order~$p^2$ 
for infinitely many~$i$.  

\proof If $G$ is elementary abelian (i.e.\ is isomorphic to a product 
of cyclic groups of order $p$) then $H^i(G)$ has exponent $p$ for $i>0$
by the K\"unneth theorem.  Conversely, if $G$ is not elementary abelian 
then $G$ contains a central subgroup of order $p$ which is not a direct 
factor, or equivalently, $C$ of order~$p$ such that the extension class 
of $G$ in $H^2(BG/C;C)$ has order $p$.  The result now follows by applying 
Corollary~2.  
\qed

The following application of Proposition 1 is new.  

\proclaim Proposition 2.  For positive integers $\alpha$,
$\beta$, $\gamma$, $\delta$ satisfying the inequalities 
$0\leq \gamma-\delta\leq {\rm min}\{\alpha,\beta\}$, let 
$G=G(\alpha,\beta,\gamma,\delta)$ be a 
$p$-group with the following presentation.   
$$G=\langle a,b,c\mid [a,c]=[b,c]=1=a^{p^\alpha}=b^{p^\beta} 
=c^{p^\gamma}, \quad [a,b]=c^{p^\delta}\rangle$$
Now let $\epsilon$ be ${\rm max}\{\alpha,\beta,2\gamma-\delta\}$.
Then for infinitely many $i$, $H^i(BG)$ has exponent
$p^\epsilon$, and at most finitely many of the groups $H^i(BG)$ 
have higher exponent.

\remark  It is easy to see that any group having a presentation 
of the above form for arbitrary $(\alpha,\beta,\gamma,\delta)$ 
also has a presentation of the above form in which the inequalities
are satisfied:  If $\gamma$ is less than $\delta$, then 
$c^{p^\delta}=c^{p^\gamma}=1$, and so in this case $G(\alpha,\beta,
\gamma,\delta)$ is isomorphic to $G(\alpha,\beta,\gamma,\gamma)$.  
On the other hand, 
the order of $[a,b]= c^{p^d}$ is bounded by the orders of $a$ and $b$
given that $c$ is central, and so the order of $c$ is bounded by
$p^{\alpha+\delta}$ and $p^{\beta+\delta}$.  Thus given a
presentation as above but not satisfying the second inequality 
we could replace $\gamma$ by $\gamma'={\rm min}\{\alpha+\delta,
\beta+\delta\}$ and obtain another presentation of the same
group.  

\proof First we recall that for any $G$ and any split surjection 
from $G$ onto $Q$, $H^*(BQ)$
occurs as a direct summand of $H^*(BG)$.  Now the above
group $G$ may be expressed as a split extension with kernel
$\langle a,c\rangle$ and quotient $\langle b\rangle\cong \Bbb
Z/p^\beta$, or as a split extension with kernel $\langle b,c \rangle$
and quotient $\langle a\rangle\cong \Bbb Z/p^\alpha$.  Hence we
deduce that $H^{2i}(BG)$ has elements of exponents $p^\alpha$ and
$p^\beta$ for all $i>0$.  
     
$G$ may also be viewed as a central extension with kernel
$\langle c\rangle$ which is isomorphic to $\Bbb Z/p^\gamma$,
and quotient isomorphic to $\Bbb Z/p^\alpha\oplus \Bbb Z/p^\beta$
generated by the images of $a$ and $b$.  The extension class of 
this extension is easily seen to have order $p^{\gamma-\delta}$,
and so it follows from Corollary 1 that for infinitely many $i$,
$H^{2i}(BG)$ contains elements of order $p^{2\gamma-\delta}$.  

For the partial converse, note that $G$ has subgroups $\langle
a,c \rangle$, $\langle b,c\rangle$, and $\langle a, b^{p^{\gamma
-\delta}}\rangle$ of index $p^\alpha$, $p^\beta$ and $p^{2\gamma
-\delta}$ respectively whose intersection is trivial, and then 
apply the following Lemma.  

\proclaim Lemma 1.  Let $G$ be a (finite) $p$-group, 
let $H_1,\ldots,H_k$ be a family of subgroups of $G$  
such that the index $|G:H_j|$ of each $H_j$ is less than or 
equal to~$p^n$, and suppose that the intersection  
$$\bigcap_{g\in G, 1\leq j\leq k} H^g_j$$ 
of the conjugates of the subgroups~$H_j$ 
is trivial.  Then $H^i(BG)$ has exponent dividing $p^n$ for all
but finitely many $i$.  

\proof Let $\Sigma_m$ be the symmetric group on $m$ symbols and
let $G_n$ be the Sylow $p$-subgroup of $\Sigma_{p^n}$.  Since 
the index of $(\Sigma_m)^p$ in $\Sigma_{mp}$ divides exactly once
by $p$ an easy induction argument using the transfer shows that 
for all~$i>0$ and all~$n$, $H^i(BG_n)$ has exponent dividing~$p^n$.  
If $H$ is a subgroup of $G$, then the kernel of the
permutation representation of $G$ on the cosets of $H$ is the
intersection of the conjugates of $H$.  Hence if $G$ has
subgroups $H_1,\ldots,H_k$ as in the statement then $G$ occurs
as a subgroup of a product of $k$ symmetric groups on at most
$p^n$ symbols, and hence as a subgroup of $(G_n)^k$.  The result
now follows from the observation due to Adem [\ade] that for any
group~$G'$ and any subgroup~$G$, the finite generation of
$H^*(BG)$ as an $H^*(BG')$-module implies that at most finitely
many of the groups $H^i(BG)$ can have 
higher exponent than the reduced cohomology~$\tilde H^*(BG')$.  
\qqed

\remark The bound given by Lemma~1 for the exponent of almost all 
of the integral cohomology groups of a $p$-group is attained 
for many groups.  For example, Proposition~2 shows that the 
bound is attained for the groups $G(\alpha,\beta,\gamma,\delta)${}.  
We were tempted to conjecture that the bound is always attained, 
but have recently found a group of order~128 whose 
index four subgroups intersect non-trivially and whose integral 
cohomology has exponent four~[\leb].  
Adem has conjectured that for~$G$ a
finite group, if $H^i(BG)$ contains elements of order~$p^n$ for
some~$i$, then it does so for infinitely many~$i$ [\ade],
and Henn has asked if this is the case [\hen].  We 
do not know if this holds for the groups 
$G(\alpha,\beta,\gamma,\delta)$.

\beginsection References.

\frenchspacing  

\def\paper#1#2/#3/#4/#5/#6/{\par\noindent [#1] #2, #3, #4 {\bf #5} #6.}

\paper \ade A. Adem/Cohomological exponents of ${\bf Z}G$-lattices/J. Pure
and Appl. Alg./58/(1989), 1--5/

\paper \add A. Adem/On the exponent of the cohomology of discrete
groups/Bull. London Math. Soc./21/(1989), 585--590/

\paper \adc A. Adem and J. F. Carlson/Discrete groups with large exponents 
in cohomology/J. Pure and Appl. Alg./66/(1990), 111--120/

\par\noindent 
[\ce] H. Cartan and S. Eilenberg, Homological Algebra, Princeton Univ. Press 
(1956).  

\paper  \evf L. Evens/The cohomology ring of a finite group/Trans. Amer. 
Math. Soc./101/(1961), 224--239/

\paper \eve L. Evens/A generalization of the transfer map in the cohomology
of groups/Trans. Amer. Math. Soc./108/(1963), 54--65/

\paper \hen H.-W. Henn/Classifying spaces with injective mod-$p$ 
cohomology/Comment. Math. Helvetici/64/(1989), 200--206/

\paper \lea I. J. Leary/A differential in the Lyndon-Hochschild-Serre 
spectral sequence/J. Pure and Appl. Alg./88/(1993), 155--168/

\par\noindent 
[\leb] I. J. Leary, Integral cohomology of some wreath products, in 
preparation.  

\paper \ven B. B. Venkov/Cohomology algebras for some classifying
spaces/Dokl. Akad. Nauk SSSR/127/(1959), 943--944 (in Russian)/

\bye